\documentclass{article}

\usepackage{arxiv}
\usepackage{matlab-prettifier}
\usepackage[utf8]{inputenc} 
\usepackage[T1]{fontenc}    
\usepackage{hyperref}       
\usepackage{url}            
\usepackage{booktabs}       
\usepackage{amsfonts}     
\usepackage{nicefrac}      
\usepackage{microtype}      
\usepackage{cleveref}      
\usepackage{lipsum}         
\usepackage{graphicx}
\usepackage[numbers]{natbib}
\usepackage{doi}

\title{Numerical computation of spectral solutions for Sturm-Liouville eigenvalue problems}

\author{Sameh Gana\\
	Department of Basic Sciences\\
	Deanship of Preparatory Year and Supporting Studies\\
	Imam Abdulrahman Bin Faisal University\\P.O. Box 1982, Dammam,34212, Saudi Arabia\\
	\texttt{sbgana@iau.edu.sa} \\
}

\begin{document}
\maketitle

\begin{abstract}
This paper focuses on the study of Sturm-Liouville eigenvalue problems. In the classical Chebyshev collocation method, the Sturm-Liouville problem is discretized to a generalized eigenvalue problem where the functions represent interpolants in suitably rescaled Chebyshev points. We are concerned with the computation of high-order eigenvalues of Sturm-Liouville problems using an effective method of discretization based on the Chebfun software algorithms with domain truncation. We solve some numerical Sturm-Liouville eigenvalue problems and demonstrate the computations' efficiency.
\end{abstract}


\section{Introduction}
The Sturm-Liouville problem arises in many applied mathematics, science, physics and engineering areas. Many biological, chemical and physical problems are described by using models based on Sturm-Liouville
equations. For example, problems with cylindrical symmetry, diffraction problems (astronomy) resolving
power of optical instruments and heavy chains. In quantum mechanics, the solutions of the radial Schrödinger
equation describe the eigenvalues of the Sturm-Liouville problem.\\
These solutions also define the bound state energies of the non-relativistic hydrogen atom.
For more applications, see ~\citep{Flugge}, ~\citep{Hanson} and ~\citep{Pryce}.\\
In this paper, we consider the Sturm-Liouville problem 
\begin{equation}
	-\frac{d}{dx}[p(x)\frac{d}{dx}]y+q(x)y=\lambda w(x)y, a\leq x\leq b
\end{equation}
\begin{equation}
	c_{a}y(a)+d_{a}y'(a)=0
\end{equation}
\begin{equation}
	c_{b}y(b)+d_{b}y'(b)=0
\end{equation}
where $p(x)>0$ , $w(x)>0$, $c_{a}$, $d_{a}$, $c_{b}$ and $d_{b}$ are constants.\\
There is a great interest in developing accurate and efficient methods of solutions for Sturm-Liouville
problems. The purpose of this paper is to determine the solution of some Sturm-Liouville problems using
the Chebfun package and to demonstrate the highest  performance of the Chebfun system compared with classical
spectral methods in solving such problems. There are many different methods for the numerical solutions
of differential equations, which include finite difference, finite element techniques, Galerkin methods, Taylor
collocation method and Chebyshef Collocation method. Spectral methods provide exponential convergence for several problems, generally with smooth solutions. The Chebfun system provides greater flexibility in solving various differential problems than the classical spectral methods. Many packages solve Sturm-Liouville problems such as MATSLISE ~\citep{Ledoux}, SLEDGE ~\citep{Fulton}, SLEIGN ~\citep{Bailey}.\\
However, these numerical methods are not suitable for the approximation of the high-index eigenvalues for
Sturm-Liouville problems.\\
The main purpose of this paper is to assert that Chebfun, along with the spectral collocation methods,
can provide accuracy, robustness and simplicity of implementation. In addition, these methods can
compute the whole set of eigenvectors and provide some details on the accuracy and numerical stability of
the results provided.\\
For more complete descriptions of the Chebyshef Collocation method and more details on the Chebfun software system, we refer to
~\citep {TrefethenF} ~\citep {TrefethenS} ~\citep {TrefethenT} ~\citep {CANUTO}  and ~\citep {FORNBERG}.

In this paper, we explain in section 2 the concept of the Chebfun System and Chebyshev Spectral Collocation methodology. Then, in section 3, some numerical examples demonstrate the method's accuracy. Finally, We end up with the conclusion section.

\section{Chebfun System and Chebyshev Spectral Collocation Methodology}

The Chebfun system, in object-oriented MATLAB, contains algorithms that amount to spectral collocation methods on Chebyshev grids of automatically determined resolution. The Chebops tools in the Chebfun system for solving differential equations are summarized in ~\citep {Hale} and ~\citep {Driscoll}. \\
The implementation of Chebops combines the numerical analysis idea of spectral collocation with the computer science idea of the associated spectral discretization matrices. The Chebfun system explained in ~\citep {TrefethenS} solves the eigenproblem by choosing
a reference eigenvalue and checks the convergence of the process.\\
The central principle of the Chebfun, along with Chebops, can accurately solve highly Sturm- Liouville problems.\\
The Spectral Collocation method for solving differential equations consists of constructing weighted interpolants of the form
~\citep{CANUTO} :

\begin{equation}
	y(x)\approx P_{N}(x)=\sum_{j=0}^{N}\frac{\alpha(x)}{\alpha(x_{j})}\phi_{j}(x)y_{j}
\end{equation}	
where $x_{j}$ for $j=0,....,N$ are interpolation nodes, $\alpha(x)$ is a weight function, $$y_{j}=y(x_{j})$$ and the interpolating functions
$\phi_{j}(x)$ satisfy $$\phi_{j}(x_{k})=\delta_{j,k}$$ and $$y(x_{k})=P_{N}(x_{k})$$ for $k= 0,....,N$.\\ Hence $P_{N}(x)$ is an interpolant of the function $y(x)$.\\ By taking $l$ derivatives of (4) and evaluating the result at the nodes  $x_{j}$, we get:
$$y^{(l)}(x_{k})\approx \sum_{j=0}^{N}\frac{d^{l}}{dx^{l}}\left[\frac{\alpha(x)}{\alpha(x_{j})}\phi_{j}(x)\right]_{x=x_{k}} \quad k= 0,....,N$$
The entries define the differentiation matrix:  
$$D^{(l)}_{k,j}= \frac{d^{l}}{dx^{l}}\left[\frac{\alpha(x)}{\alpha(x_{j})}\phi_{j}(x)\right]_{x=x_{k}}$$
The derivatives values $y^{(l)}$ are approximated  at the nodes $x_{k }$ by $ D^{(l)}y $.

The derivatives are converted to a differentiation matrix form and the
differential equation problem is transformed into a matrix eigenvalue problem.

Our interest is to compute the solutions of Sturm-Liouville problems defined in (1) with high accuracy.

First, we rewrite (1) in the following form:
\begin{equation}-\frac{d^{2}y}{dx^{2}}-\widetilde{p}(x)\frac{dy}{dx}+\widetilde{q}(x)y=\lambda \widetilde{w}(x)y,
\end{equation}
where $\widetilde{p}(x)=\frac{p'(x)}{p(x)}$, $\widetilde{q}(x)=\frac{q(x)}{p(x)}$, $\widetilde{w}(x)=\frac{w(x)}{p(x)}$ defined in the canonical interval $[-1,1]$ . Since  the differential equation is posed on $[a,b]$, it should be converted to $[-1,1]$ through the change of variable $x$ to $\frac{1}{2}((b-a)x+b+a)$.\\
The eigenfunctions 	$y(x)$ of the eigenvalue problem approximate finite terms of Chebyshev polynomials as
\begin{equation}P_{N}(x)= \sum_{j=0}^{N}\phi_{j}(x)y_{j},
\end{equation}
where the weight function $\widetilde{w}(x)$=1, $\phi_{j}(x)$ is the Chebyshev polynomial of degree $\leq N$ and $y_{j}=y(x_{j})$ , the Chebyshev collocation points are defined by: 
\begin{equation}x_{j}=cos(\frac{j\pi}{N}),\quad j=0,....,N \end{equation}
A spectral differentiation matrix for the Chebyshev collocation points is created by interpolating a polynomial through the collocation points, i.e., the polynomial
$$P_{N}(x_{k})= \sum_{j=0}^{N}\phi_{j}(x_{k})y_{j}$$
The derivatives values of the interpolating polynomial (6) at the Chebyshev collocation points (7) are:
$$ P_{N}^{(l)}(x)=\sum_{j=0}^{N}\phi_{j}^{(l)}(x_{k})y_{j}$$The differentiation matrix $D^{(l)}$ with entries 
$$D^{(l)}_{k,j}=\phi_{j}^{(l)}(x_{k})$$ is explicitly determined in ~\citep {TrefethenF} and  ~\citep {Weiden}.\\
If we rewrite equation (5) using the differentiation matrix form, we get
$$(-D_{N}^{(2)}-\widetilde{p}D_{N}^{(1)}+\widetilde{q})y=\lambda \widetilde{w}y$$
where $\widetilde{p}=diag(\widetilde{p})$, $\widetilde{q}=diag(\widetilde{q})$ and $\widetilde{w}=diag(\widetilde{w})$.\\
The boundary conditions (2) and (3) can be determined by:
$$c_{1}P_{N}(1)+d_{1}P'_{N}(1)=0$$
$$c_{-1}P_{N}(-1)+d_{-1}P'_{N}(-1)=0$$
Then the Sturm-Liouville eigenvalue problem, defined as a block operator, is transformed into a discretization
matrix diagram:
\begin{equation} 
\left( \begin{array}{c}
	-D_{N}^{(2)}-\widetilde{p}D_{N}^{(1)}+\widetilde{q} \\
	c_{1}I+d_{1}D_{N}^{(1)}\\
	c_{-1}I+d_{-1}D_{N}^{(1)}\end{array} \right)
	y
= \lambda \left( \begin{array}{c}
\widetilde{w} I\\
0	\\
0	\end{array} \right) y
\end{equation}
The approximate solutions of the Sturm-Liouville problem are defined in (1) with boundary conditions (2) and (3)
are determined by solving the generalized eigenvalue problem (8).
For more details on convergence rates, the collocation differentiation matrices and the efficiency of the
Chebyshev collocation method, see ~\cite{Weiden}.

\section{Numerical computations}
In this section, we apply the Chebyshev Spectral Collocation Methodology outlined in the previous
section, Chebfun and Chebop system described in ~\citep{TrefethenS}  and ~\citep{TrefethenT} to some Sturm-Liouville problems. We examine the accuracy and efficiency of this methodology in a selected variety of examples. In each example,
the relative error measures the technique's efficiency.
$$E_{n}=\frac{|\lambda_{n}^{exact}-\lambda_{n}|}{|\lambda_{n}^{exact}|}$$ where $\lambda_{n}^{exact}$ for $n=0,1,2,.....$ are the exact eigenvalues and $\lambda_{n}$ are the numerical eigenvalues.

\subsection {Example 1}

We consider the  Sturm-Liouville eigenvalue problem studied in ~\cite {Orszag}
\begin{equation}-\frac{d^{2}y}{dx^{2}}= \lambda w(x)y 
\end{equation}
where $w(x)>0$ , $y(0)=y(\pi)=0$.

The eigenvalue problem has an infinite number of non-trivial solutions: the eigenvalues $\lambda_{1}, \lambda_{2}, \lambda_{3},.....$ are discrete, positive real numbers and non-degenerate.
The eigenfunctions $y_{n}(x)$
associated with different eigenvalues  $\lambda_{n}$  are orthogonal with respect to the weight function $w(x)$.\\
Using the WKB theory, we approximate $\lambda_{n}$ and $y_{n}(x)$ when $n$ is large by the formulas:
$$\lambda_{n}\sim\left[ \frac{n\pi}{\int_{0}^{\pi}\sqrt{w(t)}dt}\right] ^{2}, \quad n\longrightarrow \infty$$
and 
$$y_{n}(x)\sim\left[ \int_{0}^{\pi} \frac{\sqrt{w(t)}}{2}dt\right] ^{-\frac{1}{2}}w^{-\frac{1}{4}}(x)sin\left[ n\pi\frac{\int_{0}^{x}\sqrt{w(t)}dt}{\int_{0}^{\pi}\sqrt{w(t)}dt}\right] , \quad n\longrightarrow \infty$$
We choose the weight function $w(x)=(x+\pi)^{4}$, then the Sturm-Liouville problem (9) is transformed to 
\begin{equation}
	-\frac{d^{2}y}{dx^{2}}= \lambda (x+\pi)^{4}y, \quad  y(0)=y(\pi)=0	
\end{equation}
Then, the approximate eigenvalues and eigenfunctions of the eigenvalue problem (10) are given by
$$\lambda_{n}\sim \frac{9n^{2}}{49\pi^{4}}, \quad n\longrightarrow \infty$$
and
$$y_{n}(x)\sim \sqrt{\frac{6}{7\pi^{3}}} \quad \frac{sin\left[\frac{n(x^{3}+3x^{2}\pi+3\pi^{2}x)}{7\pi^{2}}\right] }{\pi+x}, \quad n\longrightarrow \infty$$

The Chebyshev Collocation approach to solving (10) consists of constructing
the $(N + 1) \times (N + 1)$  second derivative matrix $D_{N}^{2}$ associated with the nodes (8),
but shifted from $[-1,1]$ to $[0,\pi]$. 

The incorporation of the boundary conditions
$y(0)=y(\pi)=0$ requires that the first and last rows of the matrix $D_{N}^{2}$ are removed,
as well as its first and last columns; see ~\citep{FORNBERG}.\\
The collocation approximation of the differential eigenvalue problem (10) is
now represented by the $(N - 1) \times (N - 1)$ matrix eigenvalue problem 

\begin{equation}-D_{N}^{(2)}y=\lambda \widetilde{w}y
\end{equation}  
where $\widetilde{w}= diag(w)=diag((x_{j}+\pi)^{4})$
and $y$ is the vector of approximate eigenfunction at the interior nodes $x_{j}$.

The convergence rate can be estimated theoretically. Fitting the regularity ellipse (defined in ~\citep{Tadmor}) for Chebyshev interpolation through the pole at $x=-\pi$ indicates a convergence rate of $O(\frac{1}{(3\sqrt{8})^{N}})\simeq O(0.17^{N})$.\\
The typical rate of convergence in polynomial interpolation (and also differentiation) is exponential, where
the decay rate determines the singularity's location concerning the interval. See ~\citep{Tadmor} and
~\citep{FORNBERG}.\\
We approximate the solutions of the Sturm-Liouville problem (10) by solving the Matrix eigenvalue problem (11) using a Chebfun code:
\begin{lstlisting}[style=Matlab-editor]
L = chebop(0,pi) ;
L.op = @(x,u)  -(pi+x)^-4*diff(u,2) ;
L.bc ='dirichlet';
N = 40; 
[V, D] = eigs(L,N ) ;
diag(D)
\end{lstlisting}
In table \ref{tab:table1}, we compute the first forty eigenvalues and the related relative error between the numerical
calculation and the exact solution.\\
We consider the WKB approximations by Bender and Orzag  ~\cite {Orszag} as the exact solutions for the calculations
of errors since there is no explicit form of eigenvalues.\\
We compute in Table \ref{tab:table2} the numerical values of some eigenvalues with high-index of the problem (10).
It is clear that the eigenvalues as N increases are approximately calculated with an accuracy better than
the low-order eigenvalues.\\
The numerical results in Tables \ref{tab:table1} and \ref{tab:table2} by Chebfun algorithms closely match the exact eigenvalues of the Sturm–Liouville problem in example 1.

\begin{figure}[hbt!]
	\centering 
	\caption{Some eigenfunctions of the Sturm Liouville problem in example 1 with n = 1, 20, 50, 100. }	
	\begin{tabular}{c c }
		\includegraphics[width=7cm]{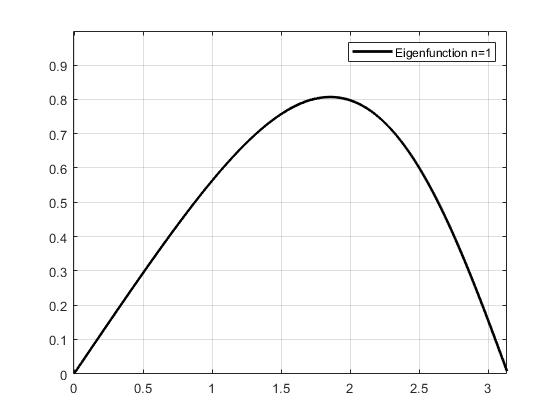} & \includegraphics[width=7cm]{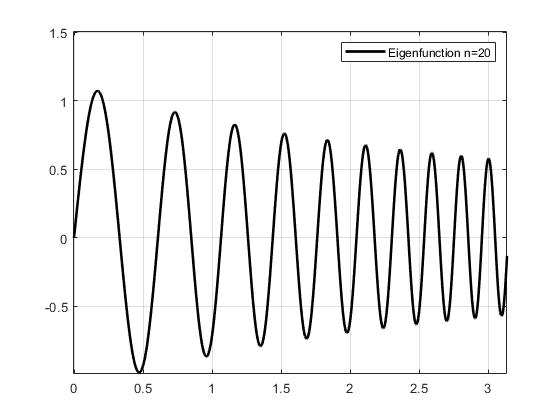}\\
		\includegraphics[width=7cm]{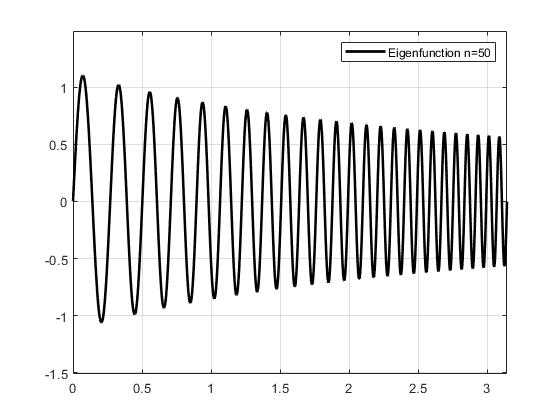}&  \includegraphics[width=7cm]{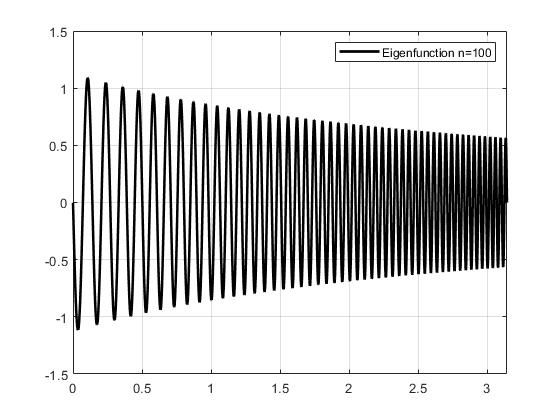}	
	\end{tabular}
\label{fig:fig1}
\end{figure}

Figure \ref{fig:fig1} shows the numerical computations of some eigenfunctions for n = 1, 20, 50 and 100.
\begin{table}[hbt!]
	\caption{Computations of the first forty eigenvalues $\lambda_{n}$ and the relative error in example 1}
	\centering\begin{tabular}{c c c c c}
		\hline
		& & & &  \\
		n\quad \quad  & $\lambda_{n}$ Current work & $\lambda_{n}^{WKB}$ \quad \quad & Relative error $E_{n}=|\frac{\lambda_{n}^{WKB}-\lambda_{n}}{\lambda_{n}^{WKB}}|$ \quad \quad& $\lambda_{n}$( ~\cite {Orszag})\quad \quad\\
		\hline
		& & & &  \\
		1	&	0.001744014		&	0.001885589		&	0.075082675	&	0.00174401 \\
		2	&	0.007348655		&	0.007542354		&	0.025681583	&	0.734865\\
		3	&	0.016752382		&	0.016970297		&	0.012840988	&	0.0167524\\
		4		& 0.029938276		&	0.030169417		&	0.00766145		& 0.0299383\\
		5	&	0.046900603		&	0.047139714		&	0.0050724 	&	0.0469006\\
		6		& 0.067636933		& 	0.067881189		&	0.003598284	 	&\\
		7		& 0.092146088			& 0.09239384		&	0.002681481		&\\
		8	&	0.120427442		&	0.120677669		&	0.002073519		&\\
		9		& 0.152480637		&	0.152732675		&	0.001650191		&\\
		10		& 0.188305458		&	0.188558858		&	0.001343874		&\\
		11		& 0.227901771		&	0.228156218		&	0.00111523		&\\
		12		& 0.271269487		&	0.271524755		&	0.00094013		&\\
		13	&	0.318408545			& 0.31866447		&	0.000803115		&\\
		14		& 0.369318906		&	0.369575361		&	0.000693919		&\\
		15		& 0.424000539		&	0.42425743		&	0.000605508		&\\
		16	&	0.482453422		&	0.482710676			&0.000532935		&\\
		17	&	0.544677542		&	0.544935099			&0.000472638		&\\
		18		&0.610672885		&	0.610930699		&	0.000422002		&\\
		19		&0.680439443		&	0.680697476		&	0.000379072		&\\
		20		&0.753977208		&	0.754235431		&	0.000342363		&0.753977\\
		21		&0.831286176		&	0.831544563			&0.00031073		&\\
		22		&0.912366343		&	0.912624871			&0.00028328		&\\
		23		&0.997217704		&	0.997476357		&	0.000259308		&\\
		24		&1.085840256		&	1.086099021		&	0.000238251		&\\
		25	&	1.178233999			&1.178492861		&	0.000219655		&\\
		26	&	1.274398929			&1.274657878		&	0.000203152		&\\
		27	&	1.374335046			&1.374594073		&	0.000188439		&\\
		28	&	1.478042348			&1.478301445		&	0.000175266		&\\
		29 & 1.585520834	&	1.585779994	&		0.000163427&\\
		30 & 1.696770504	&	1.69702972	&		0.000152747&\\
		31 & 1.811791355	&	1.812050623	&		0.00014308&\\
		32 & 1.930583389	&	1.930842703	&		0.000134301&\\
		33 & 2.053146604	&	2.053405961	&		0.000126306&\\
		34 & 2.179480999	&	2.179740395	&		0.000119003&\\
		35 & 2.309586575	&	2.309846007	&		0.000112316&\\
		36 &  2.443463331	&	2.443722796	&		0.000106176&\\
		37 & 2.581111267	&	2.581370762	&		0.000100526&\\
		38 &2.722530382	&	2.722789906	&		9.53152E-05&\\
		39 & 2.867720677&		2.867980226	&		9.0499E-05&\\
		40 & 3.016682151&		3.016941724	&		8.60386E-05& 3.01668\\
		& & & &\\
		\hline
	\end{tabular}
\label{tab:table1}
\end{table}

\begin{table}[hbt!]
	\caption{Computations of the high index eigenvalues $\lambda_{n}$ and the relative error in example 1}
	\begin{tabular}{c c c c c c c c  }
		\hline
		& & & & & & &   \\
		$N=500$ & & & & $N=1000$ & & &  \\
		\hline
		& & & & & & &   \\
		n\quad \quad  & $\lambda_{n}$  & $\lambda_{n}^{WKB}$ \quad \quad & $E_{n}$  & n\quad \quad  & $\lambda_{n}$  & $\lambda_{n}^{WKB}$ &$E_{n}$ \\
		\hline
		& & & & & & &   \\
		100	&18.56689897	&18.85588577&	0.01532608&500	&470.55992&	471.3971443	&	0.001776049\\
		150	&41.93487514	&	42.42574299	&0.011570047		&600	&	677.1850808	&678.8118879	&	0.002396551\\
		200	&75.01170433	&	75.4235431	&0.005460348	&	700		&921.8333917	&	923.9384029	&	0.002278303\\
		250	&117.0623718&	117.8492861	&	0.006677293		&750		&1059.838653	&	1060.643575	&0.0007589\\
		300	&	169.030728	&	169.702972	&0.003961298	&	800	&	1204.898974		&1206.77669		&0.001555976\\
		350	&	230.1328596	&	230.9846007	&	0.003687437	&850	&	1359.411329	&1362.337747	&	0.002148085\\
		400	&	300.5074522	&	301.6941724	&	0.00393352	&900	&	1525.004757	&1527.326748	&	0.001520297\\
		450	&	380.7355745	&381.8316869	&	0.002870669	&	950		&1700.149371	&	1701.743691	&0.000936875\\
		& & & & & & &   \\
		\hline
	\end{tabular}
\label{tab:table2}
\end{table}

\subsection {Example 2:}  
We consider the  Sturm-Liouville eigenvalue problem 
\begin{equation}-\frac{d^{2}y}{dx^{2}}+x^{4}y= \lambda y 
\end{equation}           
with the homogeneous boundary conditions $y(\pm\infty) = 0 $.

In the study of quantum mechanics, if the potential well $V(x)$
rises monotonically as $x\longrightarrow \pm \infty$, the differential equation $$\frac{d^{2}y}{dx^{2}}=(V(x)-E)y,$$
describes a particle of energy E confined to a potential well $V(x)$. 

The eigenvalue $E$ satisfies $$\int_{A}^{B}\sqrt{E-V(x)}dx=(n+\frac{1}{2})\pi,$$
where the turning points $A$ and $B$ are the two solutions to the equation
$V(x) - E = 0$. \\
The WKB eigenfunctions $y_{n}(x)$ satisfy the formula 
$$y_{n}(x)=2\sqrt{\pi}C (\frac{3}{2}S_{0})^{\frac{1}{6}}[V(x)-E]^{-\frac{1}{4}}A_{i}(\frac{3}{2}S_{0})^{\frac{2}{3}},$$
where $S_{0}=\int_{0}^{x}\sqrt{(V(t)-E)}  dt$ and $A_{i}$ is the Airy function.
( see ~\cite {Orszag})\\
Thus the eigenvalues of the problem (13) satisfy $$\lambda_{n}\sim\left[ \frac{3\gamma(\frac{3}{4})(n+\frac{1}{2})\sqrt{\pi}}{\gamma(\frac{1}{4})}\right] ^{\frac{4}{3}},\quad n\longrightarrow\infty$$
where $\gamma$ is the gamma function.

We use the Chebyshev Collocation Method by discretizing (12) in the interval $[-d,d]$ with the boundary conditions $y(-d) = y(d)=0 $:
\begin{equation}-\frac{d^{2}y}{dx^{2}}+x^{4}y= \lambda y 
\end{equation}  
with the boundary conditions $y(-d) = y(d)=0 $.\\
The collocation approximation to the differential eigenvalue problem (13) is
represented by the $(N - 1) \times (N - 1)$ matrix eigenvalue problem:
\begin{equation}-D_{N}^{(2)}y+\widetilde{q}y=\lambda y
\end{equation}
where $\widetilde{q}=diag(x_{j}^{4})$.

\begin{figure}[hbt!]
	\centering
	\caption{Relative errors E(n) of some high index eigenvalues of the Sturm-Liouville problem in example 2 with $d=10$ }
	\includegraphics[width=10cm]{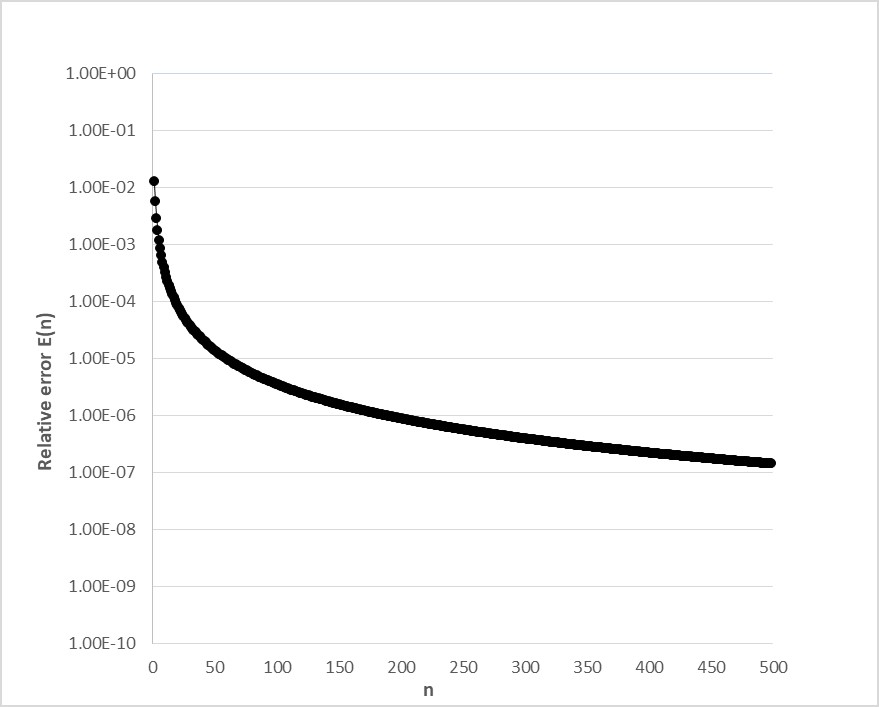}	
	\label{fig:fig2}
\end{figure}

\begin{figure}[hbt!]
	\caption{Some eigenfunctions of the Sturm Liouville problem for Example 2 with n = 10, 30, 50, 100 and $d=10$}	
	\begin{tabular}{c c}
		\includegraphics[width=9cm]{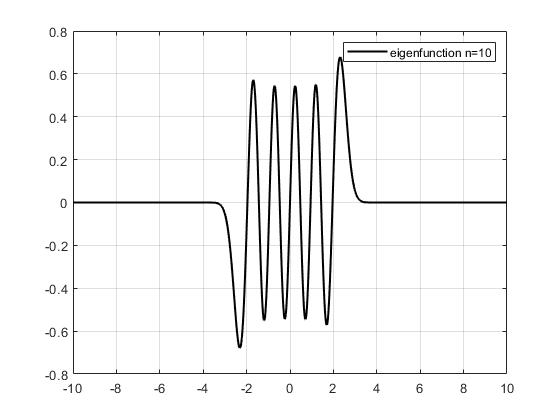} & \includegraphics[width=9cm]{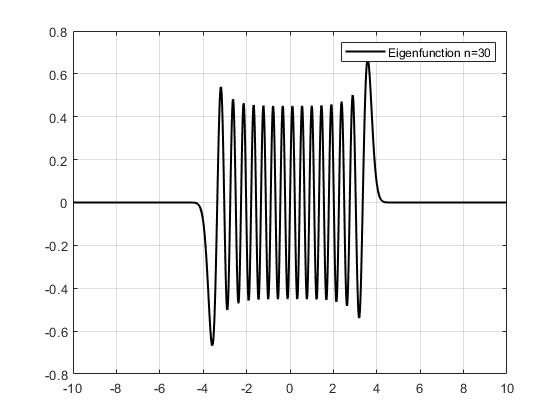}\\
		\includegraphics[width=9cm]{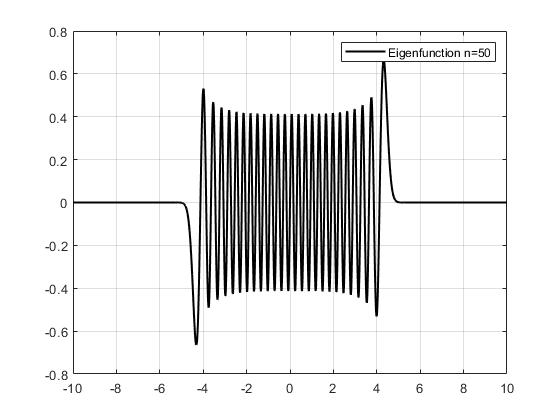} &  \includegraphics[width=9cm]{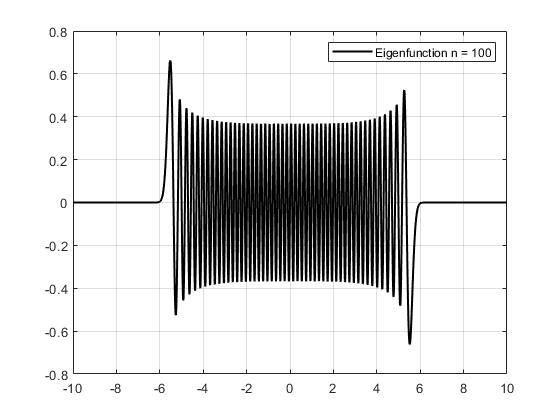}	\\
	\end{tabular}
\label{fig:fig3}
\end{figure}

We approximate the solutions of the Sturm-Liouville problem (12) by solving the Matrix eigenvalue problems
(13) and (14) using a Chebfun code:

\begin{lstlisting}[style=Matlab-editor]
d = 10;
L = chebop(-d,d);
L . op = @(x,u) - diff(u,2)+(x^4)*u;
L . bc ='dirichlet';
N = 100; 
[V, D] = eigs(L,100);
diag(D)
\end{lstlisting}

In table \ref{tab:table3}, we list the numerical results of this method by computing the first thirty eigenvalues and
the relative error between this technique and the WKB approximations by Bender and Orzag ~\cite {Orszag}. The
numerical results in Table \ref{tab:table3} show the high performance of the current technique. Figure \ref{fig:fig2} shows the numerical
computations of relative errors. These results illustrate the high accuracy and efficiency of the
algorithms.
Figure \ref{fig:fig3} shows some eigenfunctions of the sturm-liouville problem in example 2 for n = 10, 30, 50 and 100.
\begin{table}[hbt!]
	\caption{Computations of the first thirty eigenvalues $\lambda_{n}$ of and the relative error in example 2 with $d=10$}
	\centering\begin{tabular}{c c c c }
		\hline
		& & &  \\
		n\quad \quad  & $\lambda_{n}$ Current work & $\lambda_{n}^{WKB}$ \quad \quad & Relative error $E_{n}=|\frac{\lambda_{n}^{WKB}-\lambda_{n}}{\lambda_{n}^{WKB}}|$ \\
		\hline
		& & &  \\
		0&	1.0603620904849	&0.8671453264848&	2.22819357E-01\\
		1&	3.7996730297979	&3.7519199235504&	1.27276454E-02\\
		2&	7.4556979379858&	7.4139882528108	&5.62580945E-03\\
		4&	11.6447455113774&	11.6115253451971&	2.86096488E-03\\
		5&	16.2618260188517&	16.2336146927052&	1.73783391E-03\\
		6&	21.2383729182367&	21.2136533590572&	1.16526648E-03\\
		7&	26.5284711836832&	26.5063355109631&	8.35108750E-04\\
		8&	32.0985977109688&	32.0784641156416&	6.27635889E-04\\
		9&	37.9230010270330&	37.9044718450677&	4.88838943E-04\\
		10&	43.9811580972898&	43.9639483585989&	3.91451162E-04\\
		11&	50.2562545166843&	50.2401523191723&	3.20504552E-04\\
		12&	56.7342140551754&	56.7190570966241&	2.67228676E-04\\
		13&	63.4030469867205&	63.3887079062501&	2.26208751E-04\\
		14&	70.2523946286162&	70.2387714452705&	1.93955319E-04\\
		15&	77.2732004819871&	77.2602101293507&	1.68137682E-04\\
		16&	84.4574662749449&	84.4450400943621&	1.47151101E-04\\
		17&	91.7980668089950&	91.7861473252516&	1.29861467E-04\\
		18&	99.2886066604955&	99.2771452225694&	1.15448907E-04\\
		19&	106.923307381733&	106.912262402219&	1.03308819E-04\\
		20&	114.696917384982&	114.686253003331&	9.29874451E-05\\
		21&	122.604639001000&	122.594324052793&	8.41388726E-05\\
		22&	130.642068748629&	130.632075959854&	7.64956746E-05\\
		23&	138.805147911395&	138.795453260716&	6.98484745E-05\\
		24&	147.090121257603&	147.080703465973&	6.40314563E-05\\
		25&	155.493502268682&	155.484342386656&	5.89119257E-05\\
		26&	164.012043622866&	164.003124693834&	5.43826775E-05\\
		27&	172.642711962846&	172.634018745858&	5.03563379E-05\\
		28&	181.38266618577&	181.374184925625&	4.67611207E-05\\
		29&	190.229238652464&	190.220956887619&	4.35376048E-05\\	
		& & & \\
		\hline
	\end{tabular}
\label{tab:table3}
\end{table}

\subsection {Example 3:}  
\begin{figure}[hbt!]
	\caption{Some eigenfunctions of the Sturm Liouville problem for Example 3 with n = 1, 30, 50, 80.	}
	\begin{tabular}{c c}
		\includegraphics[width=9cm]{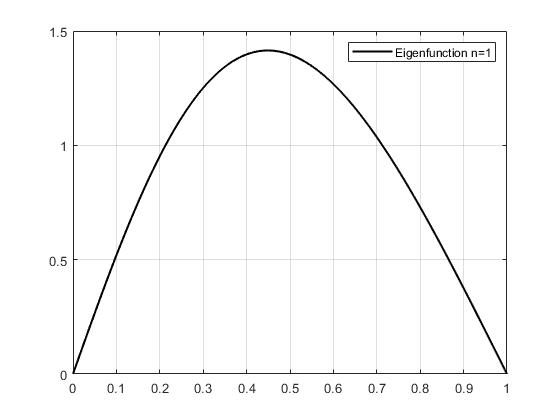} & \includegraphics[width=9cm]{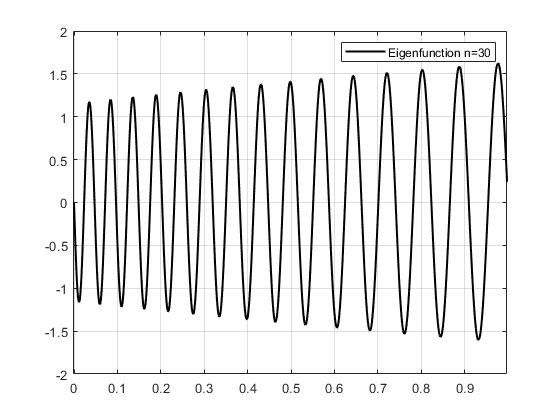}\\
		\includegraphics[width=9cm]{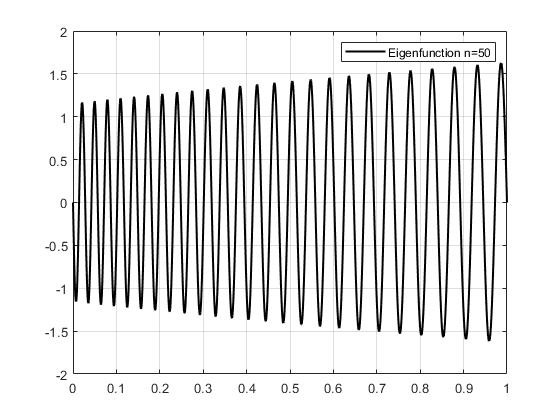} &  \includegraphics[width=9cm]{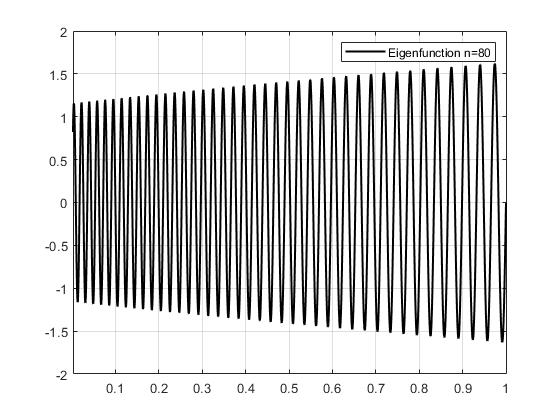}	\\
	\end{tabular}
\label{fig:fig4}
\end{figure}

We consider the  Sturm-Liouville eigenvalue problem  
\begin{equation}-\frac{d^{2}y}{dx^{2}}= \frac{\lambda}{(1+x)^{2}}y 
\end{equation}  
with the boundary conditions $y(0)=y(1)=0$.
The exact eigenvalues of the problem (15) satisfy the explicit formula (see ~\cite {Akulenko}):
$$\lambda_{n}=\frac{1}{4}+(\frac{\pi n}{ln2})^{2}, \quad n=1,2,3,...$$
The eigenfunctions $y_{n}(x)$
associated with different eigenvalues  $\lambda_{n}$ are given by:
$$y_{n}(x)=\sqrt{1+x}\quad sin(\frac{\pi n}{ln2} ln(1+x)).$$
The collocation approximation of the differential eigenvalue problem (15) is
now represented by the $(N - 1) \times (N - 1)$ matrix eigenvalue problem 
\begin{equation}-D_{N}^{(2)}y=\lambda \widetilde{w}y
\end{equation}  
where $\widetilde{w}=diag(\frac{1}{(1+x_{j})^{2}})$.\\
Now, the approximate eigenvalues of the Sturm-Liouville problem (15) are obtained by solving the matrix
eigenvalue problem (16) using the Chebyshev Spectral Collocation technique based on Chebfun and Chebop
codes.\\
In table \ref{tab:tab4}, we compute the first thirty eigenvalues and the related absolute error between the numerical
calculation and the exact solution
$$E_{n}=\frac{|\lambda_{n}^{exact}-\lambda_{n}|}{|\lambda_{n}^{exact}|}.$$ 
The eigenvalues obtained are extremely close to the exact eigenvalues. The results
show significant improvement in the convergence.\\
In figure \ref{fig:fig4}, we plot some eigenfunctions in example 3 for $ n=1, n=30 , n=50 $ and $n=80$.
\begin{table}[hbt!]
	\caption{Computations of the first thirty eigenvalues $\lambda_{n}$ and the relative error in example 3}
	\centering\begin{tabular}{c c c c }
		\hline
		& & &  \\
		n\quad \quad  & $\lambda_{n}$ Current work & $\lambda_{n}^{exact}$ \quad \quad & Relative error $E_{n}=\frac{|\lambda_{n}^{exact}-\lambda_{n}|}{|\lambda_{n}^{exact}|}$ \\
		\hline
		& & &  \\
		1&	20.79228845522		&20.79228845517	&	2.40469E-12\\
		2&	82.4191538209		&82.41915382087		&3.63913E-13\\
		3	&185.13059609701	&	185.13059609709	&	4.32157E-13\\
		4&	328.92661528358	&	328.92661528388		&9.11961E-13\\
		5	&	513.8072113806	&	513.80721137895		&3.21132E-12\\
		6	&	739.77238438806	&	739.77238439069		&3.55512E-12\\
		7		&1006.82213430597	&	1006.82213430587	&	9.9243E-14\\
		8	&	1314.95646113432	&	1314.95646113354	&	5.93335E-13\\
		9	&	1664.17536487313	&	1664.17536487301	&	7.20502E-14\\
		10		&2054.47884552238	&	2054.47884552193	&	2.18994E-13\\
		11		&2485.86690308208	&	2485.86690308175	&	1.32756E-13\\
		12		&2958.33953755223	&	2958.33953755204	&	6.41739E-14\\
		13		&3471.89674893283	&	3471.89674893313	&	8.6339E-14\\
		14		&4026.53853722387	&	4026.53853722418	&	7.70655E-14\\
		15	&	4622.26490242536	&	4622.26490242571	&	7.578E-14\\
		16		&5259.0758445373	&	5259.07584453724	&	1.13997E-14\\
		17	&	5936.97136355969	&	5936.97136355928	&	6.90036E-14\\
		18	&	6655.95145949252	&	6655.95145949275	&	3.46947E-14\\
		19	&	7416.0161323358	&	7416.01613233622	&	5.65889E-14\\
		20	&	8217.16538208953	&	8217.16538208989	&	4.39108E-14\\
		21	&	9059.39920875371	&	9059.3992087539		&2.09559E-14\\
		22	&	9942.71761232833	&	9942.71761232853	&	2.00991E-14\\
		23	&	10867.1205928134	&	10867.1205928132	&	1.83894E-14\\
		24	&	11832.6081502089	&	11832.6081502087	&	1.68889E-14\\
		25	&	12839.1802845148	&	12839.1802845162	&	1.09127E-13\\
		26	&	13886.8369957313	&	13886.8369957319	&	4.31718E-14\\
		27	&	14975.5782838581	&	14975.5782838589	&	5.33776E-14\\
		28	&	16105.4041488954	&	16105.4041488943	&	6.82455E-14\\
		29	&	17276.3145908432	&	17276.3145908419	&	7.53159E-14\\
		30	&	18488.3096097014	&	18488.3096097008	&	3.25471E-14\\		
		& & & \\
		\hline
	\end{tabular}
\label{tab:tab4}
\end{table}
\section*{Conclusion}
The numerical computations prove the efficiency of the technique based on the Chebfun and Chebop systems.
This technique is unbeatable regarding the accuracy, computation speed, and information it provides on
the accuracy of the computational process. Chebfun provides greater flexibility compared to classical
spectral methods.
However, in the presence of various singularities, the maximum order of approximation can be reached
and the Chebfun issues a message that warns about the possible inaccuracy of the computations.
The methodology can be used to obtain high-accurate solutions to other Sturm-Liouville problems, generalized differential equations involving higher order derivatives and non-linear partial differential equations in
multiple space dimensions.

\bibliographystyle{unsrtnat}
\clearpage
\bibliography{references}

\end{document}